\documentclass[11pt]{book}
\usepackage{harvard}
\usepackage{amsthm}
\usepackage{amsfonts}
\usepackage{mhequ}
\usepackage{mathrsfs}

\newtheorem{theorem}{Theorem}[section]
\newtheorem{corollary}[theorem]{Corollary}
\newtheorem{lemma}[theorem]{Lemma}
\newtheorem{remark}[theorem]{Remark}
\newtheorem{proposition}[theorem]{Proposition}
\newtheorem{definition}[theorem]{Definition}
\newtheorem{example}[theorem]{Example}

\makeatletter
\let\@old@mean\mean\def\mean{\mathfrak{m}}
\let\@old@d\d\def\d{\partial}
\let\@old@cM\cM\def\cM{\mathscr{M}}
\let\@old@cP\cP\def\cP{\mathscr{P}}
\let\@old@cF\cF\def\cF{\mathscr{F}}
\let\@old@P\P\def\P{\mathbf{P}}
\let\@old@E\E\def\E{\mathbf{E}}
\let\@old@R\R\def\R{\mathbf{R}}
\let\@old@Z\Z\def\Z{\mathbf{Z}}
\let\@old@N\N\def\N{\mathbf{N}}
\let\@old@CC\CC\def\CC{\mathcal{C}}
\let\@old@CH\CH\def\CH{\mathcal{H}}
\let\@old@CB\CB\def\CB{\mathcal{B}}
\let\@old@CF\CF\def\CF{\mathcal{F}}
\let\@old@CM\CM\def\CM{\mathcal{M}}
\let\@old@CP\CP\def\CP{\mathcal{P}}
\let\@old@CQ\CQ\def\CQ{\mathcal{Q}}
\let\@old@CS\CS\def\CS{\mathcal{S}}
\let\@old@CX\CX\def\CX{\mathcal{X}}
\let\@old@CY\CY\def\CY{\mathcal{Y}}
\let\@old@CW\CW\def\CW{\mathcal{W}}
\let\@old@TV\TV\def\TV{\mathrm{TV}}
\let\@old@eref\eref\def\eref#1{(\ref{#1})}
\let\@old@scal\scal\def\scal#1{\langle{#1}\rangle}
\makeatother

\begin{document}

\title{Ergodic properties of a class of non-Markovian processes}
\author{M. Hairer}

\maketitle

\chapter[non-Markovian processes]{Ergodic properties of a class of non-Markovian processes}

\subsection*{Abstract}

We study a fairly general class of time-homogeneous stochastic evolutions driven by noises
that are not white in time. As a consequence, the resulting processes do not have
the Markov property. In this setting, we obtain constructive criteria for the uniqueness of
stationary solutions that are very close in spirit to the existing criteria for Markov processes.

In the case of discrete time, where the driving noise consists of a stationary sequence
of Gaussian random variables, we give optimal conditions on the spectral measure 
for our criteria to be applicable. In particular, we show that under a certain assumption on the spectral
density, our assumptions can be checked in virtually the same way as one would check
that the Markov process obtained by replacing the
driving sequence by a sequence of independent identically distributed
Gaussian random variables is strong Feller and topologically irreducible.
The results of the present article
extend those obtained previously in the continuous time context of 
 diffusions driven by fractional Brownian 
motion. 

\section{Introduction}

Stochastic processes have been used as a powerful modelling tool for decades in situations
where the evolution of a system has some random component, be it intrinsic or to model the interaction
with a complex environment. In its
most general form, a stochastic process describes the evolution $X(t,\omega)$ 
of a system, where $t$ denotes the time parameter and 
$\omega$ takes values in some probability space and abstracts the
`element of chance' describing the randomness of the process.

In many situations of interest, the evolution of the system can be described (at least informally) 
by the solutions of an evolution equation of the type
\minilab{e:main}
\begin{equ}[e:maincont]
{dx \over dt} = F(x, \xi)\;,
\end{equ}
where $\xi$ is the `noise' responsible for the randomness in the evolution.
In the present article, we will not be interested in the technical subtleties arising
from the fact that the time parameter $t$ in \eref{e:maincont} takes continuous values.
We will therefore consider its discrete analogue
\minilab{e:main}
\begin{equ}[e:maindiscr]
x_{n+1} = F(x, \xi_n)\;,
\end{equ}
were $\xi_n$ describes the noise acting on the system between times $n$ and $n+1$.
Note that \eref{e:maincont} can always be reduced to \eref{e:maindiscr} by allowing
$x_n$ to represent not just the state of the system at time $n$, but its evolution over the whole
time interval $[n-1, n]$.
We were intentionally vague about the precise meaning of the symbol $x$ in the right hand side of
\eref{e:main} in order to suggest that there are situations where it makes sense to let the
right hand-side depend not only on the current state of the system, but on the whole collection
of its past states as well.

The process $x_n$ defined by a recursion of the type \eref{e:maindiscr} has the Markov property
if both of the following properties hold:
\begin{itemize}
\item[\textit{a.}] The noises $\{\xi_n\}_{n \in \Z}$ are mutually independent.
\item[\textit{b.}] For a fixed value of $\xi$, the  function $x \mapsto F(x,\xi)$ depends only on the last state of the system.
\end{itemize}
In this article, we will be interested in the study of recursion relations of the type
\eref{e:maindiscr} when condition \textit{b.}\ still holds, but the Markov property
is lost because condition \textit{a.}\ fails to hold. Our main focus will be on the ergodic
properties of \eref{e:maindiscr}, with the aim of providing concrete conditions that ensure the
uniqueness (in law) of a stationary sequence of random variables $x_n$ satisfying a given
recursion of the type \eref{e:maindiscr}.

Many such criteria exist for Markov processes and we refer to \cite{MT93} for 
a comprehensive overview of the techniques developed in this regard over the past
seven decades. The aim of the present article will be to present a framework in which 
recursions of the type \eref{e:maindiscr} can be studied and such that several existing 
ergodicity results for Markov processes have natural equivalents whose assumptions can
also be checked in similar ways. This framework (which should be considered
as nothing but a different way of looking at random dynamical systems, together with some
slightly more restrictive topological assumptions) was developed in
\cite{Hai05} and further studied in \cite{HaiOha07} in order to treat the ergodicity of stochastic
differential equations driven by fractional Brownian motion. The main novelty of the present article
is to relax a number of assumptions from the previous works and to include a detailed
study of the discrete-time case when the driving noise is Gaussian.

The remainder of this article is organised as follows. After introducing our notations at the end of this
section, we will introduce in Section~\ref{sec:framework} the framework studied in the present
article. We then proceed in  Section~\ref{sec:RDS} to a comparison of this framework with that of 
random dynamical systems. In Section~\ref{sec:ergodicity} we recall a few general ergodicity criteria
for Markov processes and give a very similar criterion that can be applied in our framework.
In Section~\ref{sec:Gauss} we finally study in detail the case of a system driven by a (time-discrete)
stationary sequence of Gaussian random variables. We derive an explicit condition on the
spectral measure of the sequence that ensures that such a system behaves qualitatively
like the same system driven by an i.i.d.\ sequence of Gaussian random variables.

\subsection{Notations}

The following notations will be used throughout this article. Unless stated otherwise, measures will always be Borel
measures over Polish (\textit{i.e.}\ metrisable, complete, separable) spaces and they will always be positive.
We write $\cM_+(\CX)$ for the set of all such measures on the space $\CX$ and $\cM_1(\CX)$ for the subset
of all probability measures. We write $\mu \approx \nu$ to indicate that $\mu$ and $\nu$ are equivalent
(\textit{i.e.}\ they are mutually absolutely continuous, that is they have the same negligible sets) 
and $\mu \perp \nu$ to indicate that they are mutually singular.

Given a map $f \colon \CX \to \CY$ and a measure $\mu$ on $\CX$, we denote by $f^* \mu$
the push-forward measure $\mu \circ f^{-1}$ on $\CY$. Given a product space $\CX \times \CY$, we will use the
notation $\Pi_\CX$ and $\Pi_\CY$ to denote the projections onto the two factors. For infinite products
like $\CX^\Z$, $\CX^{\Z_-}$ or $\CX^\N$, we denote by $\Pi_n$ the projection onto the $n$th factor
($n$ can be negative in the first two cases).

We will also make use of the concatenation operator $\sqcup$ from $\CX^{\Z_-} \times \CX^n$
to $\CX^{\Z_-}$ defined in the natural way by
\begin{equ}
(w \sqcup w')_k = \left\{\begin{array}{rl} w_{k+n} & \text{if $k \le -n$,} \\ w'_{k+n} & \text{otherwise.} \end{array}\right. 
\end{equ}

Finally, given a Markov transition probability $\CP\colon \CX \to \cM_1(\CX)$, we will use the same
symbol for the associated Markov operator acting on observables $\phi \colon \CX \to \R$ by
$\bigl(\CP \phi\bigr)(x) = \int_\CX \phi(y) \,\CP(x,dy)$, and the dual operator acting on probability
measures $\mu$ by $\bigl(\CP \mu\bigr)(A) = \int_\CX \CP(x,A)\, \mu(dx)$.

\subsection*{Acknowledgements}
{\small
The author would like to thank J.~Mattingly for countless conversations on the topic of this article,
as well as S.~Assing and J.~Vo\ss\ for their helpful comments that hopefully lead to greatly improved clarity
and the correction of many misprints. The report of a careful referee lead to further clarity.
This work was supported by EPSRC fellowship EP/D071593/1.
}

\section{Skew-products}
\label{sec:framework}

Whatever stochastic process $X$ one may wish to consider, it is always possible to
turn it into a Markov process by adding sufficiently many `hidden' degrees of freedom to the state space.
For example, one can take the state space large enough to contain all possible information
about the past of $X$, as well as all possible information on the future of
the driving noise $\xi$.
The evolution \eref{e:main} is then deterministic, with all randomness injected once and for
all by drawing $\xi$ initially according to the appropriate distribution.
This is the point of view of random dynamical systems explained in more detail in
Section~\ref{sec:RDS} below.

On the other hand, one could take a somewhat smaller `noise space' that contains
only information about the \textit{past} of the driving noise $\xi$. In this case, the evolution
is no longer deterministic, but it becomes a skew-product between a Markovian
evolution for the noise (with the transition probabilities given informally by the conditional
distribution of the `future' given the `past') and a deterministic map that solves \eref{e:main}.
This is the viewpoint that was developed in \cite{Hai05,HaiOha07} and will be studied 
further in this article.

The framework that will be considered here is the following. Let $\CW$ and $\CX$ be two
Polish spaces that will be called the `noise space $\CW$'
and the `state space $\CX$' respectively, let $\CP$ be a Markov transition kernel on $\CW$, and let
$\Phi\colon \CW \times \CX \to \CX$ be an `evolution map'. Throughout this article,
we will make the following standing assumptions:
\begin{itemize}
\item[1.] There exists a probability measure $\P$ on $\CW$ which is invariant for $\CP$ and such that
the law of the corresponding stationary process is ergodic under the shift map.
\item[2.] The map $\Phi \colon \CW \times \CX \to \CX$ is continuous in the product topology.
\end{itemize}
We will also occasionally impose some regularity of the kernel $\CP(w,\cdot\,)$ as a function of $w$.
We therefore state the following property which will not always be assumed to hold:
\begin{itemize}
\item[3.] The transition kernel $\CP$ is Feller, that is the function $\CP \phi$ defined
by $\bigl(\CP \phi\bigr)(w) = \int_\CW \phi(w')\CP(w,dw')$ is continuous as soon as $\phi$ is continuous.
\end{itemize}

There are two objects that come with a construction such as the one above. First, we can
define a Markov transition operator $\CQ$ on $\CX \times \CW$ by
\begin{equ}[e:evolSkew]
\bigl(\CQ \phi\bigr)(x,w) = \int_{\CW} \phi(\Phi(x,w'), w')\CP(w, dw')\;.
\end{equ}
In words, we first draw an element $w'$ from the noise space according to the law $\CP(w,\cdot\,)$
and we then update the state of the system with that noise according to $\Phi$. We also introduce
a `solution map' $\CS\colon \CX\times \CW \to \cM_1(\CX^{\N})$
 that takes as arguments an initial condition 
$x\in \CX$ and an `initial noise' $w$ and returns the law of the corresponding solution process,
that is the marginal on $\CX$ of the law of the Markov process starting at $(x,w)$ with transition
probabilities $\CQ$.

The point of view that we take in this article is that $\CS$ encodes all the `physically relevant' 
part of the evolution \eref{e:main}, and that the particular choice of noise space is just a mathematical
tool. This motivates the introduction of an equivalence relation between probability measures
on $\CX \times \CW$ by
\begin{equ}[e:equivalence]
\mu \simeq \nu \quad\Leftrightarrow\quad \CS \mu = \CS \nu\;.
\end{equ}
(Here, we used the shorthand $\CS\mu = \int \CS(x,w)\, \mu(dx,dw)$.)
In the remainder of this article, when we will be looking for criteria that ensure the uniqueness of
the invariant measure for $\CQ$, this will always be understood to hold
up to the equivalence relation \eref{e:equivalence}.

\begin{remark}
The word `skew-product' is sometimes used in a slightly different way in the literature.
In our framework, given a realisation of the noise, that is a realisation of a Markov process on
$\CW$ with transition probabilities $\CP$, the evolution in $\CX$ is purely deterministic. This
is different from, for example, the skew-product decomposition of Brownian motion where, given
one realisation of the evolution of the radial part, the evolution of the angular part is still random. 
\end{remark}

\subsection{Admissible measures}

We consider the invariant measure $\P$ for the noise process to be fixed. Therefore, we
will usually consider measures on $\CX \times \CW$ such that their projections on $\CW$
are equal to $\P$. Let us call such probability measures \textit{admissible} and let us denote
the set of admissible probability measures by $\cM_\P(\CX)$. Obviously, the Markov
operator $\CQ$ maps the set of admissible probability measures into itself. Since we assumed
that $\CW$ is a Polish space, it is natural to endow $\cM_\P(\CX)$ with the topology
of weak convergence. This topology is preserved by $\CQ$ if we assume that $\CP$ is Feller:

\begin{lemma}
If $\Phi$ is continuous and $\CP$ is Feller, then the Markov transition operator $\CQ$ is also Feller
and therefore continuous in the topology of weak convergence on $\CX \times\CW$.\qed
\end{lemma}

The proof of this result is straightforward and of no particular interest, so we leave it as an exercise.

There are however cases of interest in which we do not wish to assume that $\CP$
is Feller. In this case, a natural topology for the space $\cM_\P(\CX)$ is given by the
`narrow topology', see \cite{Val90Young,Crau02RPM}. In order to define this topology,
denote by $\CC_\P(\CX)$ the set of functions $\phi \colon \CX \times \CW \to \R$ such that
$x \mapsto \phi(x,w)$ is bounded and continuous for every $w \in \CW$, $w \mapsto \phi(x,w)$
is measurable for every $x \in \CX$, and $\int_\CW \sup_{x \in \CX} |\phi(x,w)|\,\P(dw) < \infty$.
The narrow topology on  $\cM_\P(\CX)$ is then the coarsest topology such that the map
$\mu \mapsto \int \phi(x,w)\, \mu(dx,dw)$ is continuous for every $\phi \in \CC_\P(\CX)$.
Using Lebesgue's dominated convergence theorem, it is straightforward to show that 
$\CQ$ is continuous in the narrow topology without requiring any assumption besides the
continuity of $\Phi$.

An admissible probability measure $\mu$ is now called an \textit{invariant measure} 
for the skew-product 
$(\CW,\P,\CP,\CX,\Phi)$ if it is an invariant measure for $\CQ$, that is if $\CQ\mu = \mu$.
We call it a \textit{stationary measure} if $\CQ \mu \simeq \mu$, that is if the law of the 
$\CX$-component of the Markov process with transition probabilities $\CQ$ starting from $\mu$
is stationary. Using the standard Krylov-Bogoliubov argument, one shows that

\begin{lemma}
Given any stationary measure  $\mu$ as defined above, there exists an invariant measure
$\hat \mu$ such that $\hat \mu \simeq \mu$.
\end{lemma}

\begin{proof}
Define a sequence of probability measures $\mu_N$ on $\CX \times \CW$ by
$\mu_N = {1\over N} \sum_{n=1}^N \CQ^n \mu$.
Since, for every $N$, the marginal of $\mu_N$ on $\CW$ is equal to $\P$ and the marginal on $\CX$
is equal to the marginal of $\mu$ on $\CX$ (by stationarity), this sequence is tight in the narrow topology \cite{Crau02RPM}. 
It therefore has
at least one accumulation point $\hat \mu$ and the continuity of $\CQ$ in the narrow topology
ensures that 
$\hat \mu$ is indeed an invariant measure for $\CQ$.
\end{proof}

The aim of this article is to present some criteria that allow to show the uniqueness
\textit{up to the equivalence relation \eref{e:equivalence}} of the invariant measure
for a given skew-product. The philosophy that we will pursue is not to apply
existing criteria to the Markov semigroup $\CQ$. This is because, in typical situations
like a random differential equation driven by some stationary Gaussian process, the noise
space $\CW$ is very `large' and so the Markov operators $\CP$ and $\CQ$
typically do not have any of the `nice' properties (strong Feller
property, $\psi$-irreducibility, etc.) that are often required in the ergodic theory of Markov 
processes.

\subsection{A simple example}
\label{sec:example}

In this section, we give a simple example that illustrates the fact that it is possible in some
situations to have non-uniqueness of the invariant measure for $\CQ$, even though
$\P$ is ergodic and 
one has uniqueness up to the equivalence relation \eref{e:equivalence}. Take 
$\CW = \{0,1\}^{\Z_-}$ and define the `concatenation and shift' map $\Theta\colon \CW \times \{0,1\} \to \CW$ by
\begin{equ}
\Theta(w,j)_n = \left\{\begin{array}{cl} w_{n+1} & \text{for $n < 0$,} \\ j & \text{for $n = 0$.} \end{array}\right.
\end{equ}
Fix $p \in (0,1)$, let $\xi$ be a random variable that takes the values $0$ and $1$ with probabilities $p$
and $1-p$ respectively, and define the transition probabilities $\CP$ by 
\begin{equ}
\CP(w,\cdot\,) \stackrel{\mbox{\tiny law}}{=} \Theta(w,\xi)\;.
\end{equ}
We then take as our state space $\CX = \{0,1\}$ and we define an evolution $\Phi$ by
\begin{equ}[e:defphiex]
\Phi(x,w) = \left\{\begin{array}{cl} w_0 & \text{if $x = w_{-1}$,} \\ 1-w_0 & \text{otherwise.} \end{array}\right.
\end{equ}It is clear that there are two extremal invariant measures for this evolution. One of them charges
the set of pairs $(w,x)$ such that $x = w_0$, the other one charges the set of pairs such that $x = 1-w_0$. 
(The projection of these measures onto $\CW$ is Bernoulli with parameter $p$ in both cases.)
However, if $p = {1\over 2}$, both invariant measures give raise to the same stationary process in $\CX$, 
which is just a sequence of independent Bernoulli random variables.

\subsection{An important special case}
\label{sec:prod}

In most of the remainder of this article, we are going to 
focus on the following particular case of the setup described above.
Suppose that there exists a Polish space $\CW_0$ which carries all the information
that is needed in order to reconstruct the dynamic of the system from one time-step to the
next. We then take $\CW$ of the form $\CW = \CW_0^{\Z_-}$ (with the product topology)
and we assume that $\Phi$ is of the form $\Phi(x, w) = \Phi_0(x, w_0)$ for some jointly
continuous function $\Phi_0 \colon \CX \times \CW_0 \to \CX$. Here, we use the notation
$w = (\ldots, w_{-1}, w_0)$ for elements of $\CW$. Concerning the transition probabilities
$\CP$, we fix a Borel measure $\P$ on $\CW$ which is invariant and ergodic for the shift 
map\footnote{Recall that a probability measure $\mu$ is ergodic for a map $T$ leaving $\mu$ invariant if 
all $T$-invariant measurable sets are of $\mu$-measure $0$ or $1$.}
$(\Theta w)_n = w_{n-1}$, and we define a measurable map $\bar \CP \colon \CW \to \cM_1(\CW_0)$
as the regular conditional probabilities of $\P$ under the splitting $\CW \approx \CW \times \CW_0$.

The transition probabilities $\CP(w,\cdot\,)$ are then constructed as the push-forward of 
$\bar \CP(w,\cdot\,)$ under the concatenation map $f_w \colon \CW_0 \to \CW$ given by  $f_w(w') = w \sqcup w'$. 
Since we assumed that $\P$ is shift-invariant, it follows from the construction that
it is automatically invariant for $\CP$.

Many natural situations fall under this setup, even if they do not look so at first sight. For example,
in the case of the example from the previous section, one would be tempted to take
$\CW_0 = \{0,1\}$. This does not work since the function $\phi$ defined
in \eref{e:defphiex} depends not only on $w_0$ but also on $w_{-1}$. However, one can
choose $\CW_0 = \{0,1\}^2$ and identify $\CW$ with the subset of all sequences $\{w_n\}_{n \le 0}$ in
$\CW_0^{\Z_-}$ such that $w_n^2 = w_{n+1}^1$ for every $n < 0$.

\section[Markov processes versus RDS]{Skew-products of Markov processes 
versus random dynamical systems}
\label{sec:RDS}

There already exists a mature theory which was developed precisely
in order to study systems like \eref{e:main}. The theory in question is of course 
that of random dynamical
systems (RDS in the sequel), which was introduced under this name in the nineties by Arnold and then 
developed further by
a number of authors, in particular Caraballo, Crauel, 
Debussche, Flandoli, Robinson, Schmalfu\ss, and many others. Actually, skew-products of flows had
been considered by authors much earlier, see for example
\cite{SS73,SS77}, or the monograph \cite{Kif86}, but previous authors usually made very restrictive
assumptions on the structure of the noise, either having independent noises at each step or some
periodicity or quasi-periodicity. We refer
to the monograph \cite{Arn98} for a thorough exposition of the theory, but for the sake
of completeness, we briefly recall the main framework here. For simplicity,
and in order to facilitate comparison with the alternative 
framework presented in this article, we restrict ourselves
to the case of discrete time. An RDS consists of a dynamical system $(\Omega,\bar\P,\Theta)$, 
(here $\bar\P$ is a probability on the measurable space $\Omega$ which is both 
invariant and ergodic for the map $\Theta\colon\Omega\to\Omega$) together with a
`state space' $\CX$ and a map $\bar\Phi\colon \Omega\times \CX \to \CX$. For every initial condition
$x \in \CX$, this allows to construct a stochastic process $X_n$ over $(\Omega,\bar\P)$,
viewed as a probability space, by
\begin{equ}
X_0(\omega) = x\;,\quad
X_{n+1}(\omega) = \bar\Phi\bigl(\Theta^n \omega, X_{n}(\omega)\bigr)\;.
\end{equ}
Note the similarity with \eref{e:evolSkew}. The main difference is that the evolution $\Theta$ on
the `noise space' $\Omega$ is \textit{deterministic}. This means that an element of
$\Omega$ must contain all possible information on the future of the noise
driving the system. In fact, one can consider an RDS as a dynamical system $\hat \Phi$
over the product space $\Omega \times \CX$ via
\begin{equ}
\hat \Phi(x, \omega) = \bigl(\Theta \omega, \bar\Phi(x, \omega)\bigr)\;.
\end{equ}
The stochastic process $X_n$ is then nothing but the projection on $\CX$ of a `typical'
orbit of $\hat \Phi$. An \textit{invariant measure} for an RDS $(\Omega,\bar\P,\Theta, \CX, \bar\Phi)$
is a probability measure $\mu$ on $\Omega \times \CX$ which is invariant under $\hat \Phi$
and such that its marginal on $\Omega$ is equal to $\bar\P$. It can be shown \cite{Arn98} that such measures can be described by their `disintegration' over $(\Omega, \bar \P)$, which is a map $\omega \mapsto \mu_\omega$ from $\Omega$ into 
the set of probability measures on $\CX$.

Consider the example of an elliptic diffusion $X_t$ on a compact manifold $\CM$. Let us be
even more concrete and take for $X_t$ a simple Brownian motion and for $\CM$ the unit circle $S^1$.
When considered from the point of view of the Markov semigroup $\CP_t$ generated by $X_t$, 
it is straightforward to show that there exists a unique invariant probability measure
for $\CP_t$. In our example, this invariant measure is of course simply the
Lebesgue measure on the circle.

Consider now $X_t$ as generated from a random dynamical system (since we focus on the discrete
time case, choose $t$ to take integer values). At this stage,
we realise that we have a huge freedom of choice when it comes to finding an underlying dynamical
system $(\Omega, \bar\P, \Theta)$ and a map $\bar\Phi \colon \Omega \times S^1 \to S^1$ such that the
corresponding stochastic process is equal to our Brownian motion
$X_t$. The most immediate choice would be to
take for $\Omega$ the space of real-valued continuous functions vanishing at the origin,
$\bar\P$ equal to the Wiener measure, $\Theta$ the shift map $\bigl(\Theta B\bigr)(s) = B(s+1)-B(1)$, and
$\bar\Phi(x, B) = x + B(1)$. In this case, it is possible to show that there exists a unique invariant
measure for this random dynamical systems and that this invariant measure is
equal to the product measure of $\bar\P$ with Lebesgue measure on $S^1$, so that $\mu_\omega$ is
equal to the Lebesgue measure for every $\omega$.

However, we could also have considered $X_t$ as
the solution to the stochastic differential equation
\begin{equ}
dX(t) = \sin(k X(t))\, dB_1(t) + \cos(k X(t))\, dB_2(t)\;,
\end{equ}
where $B_1$ and $B_2$ are two independent Wiener processes and $k$ is an arbitrary integer.
In this case, it turns out \cite{Lej87,Crau02} that there are \textit{two} invariant measures $\mu^+$ and $\mu^-$ for the corresponding
random dynamical system. Both of them are such that, for almost every $\omega$, 
$\mu_\omega$ is equal to a sum of $k$
$\delta$-measures of weights $1/k$. Furthermore, the map $\omega \mapsto \mu^-_\omega$ is measurable with respect to
the filtration generated by the increments of $B_i(t)$ for negative $t$, whereas 
 $\omega \mapsto \mu^+_\omega$ is measurable with respect to
the filtration generated by the increments of $B_i(t)$ for positive $t$.

What this example makes clear is that while the ergodic theory of $X_t$ considered as
a Markov process focuses on the long-time behaviour of \textit{one instance} of $X_t$ started at an
arbitrary but fixed initial condition, the theory of random dynamical systems instead focuses
on the (potentially much richer) simultaneous long-time behaviour of \textit{several instances} of 
$X_t$ driven by the same instance of the noise. Furthermore, it shows that  a random dynamical
system may have invariant measures that are `unphysical' in the sense that they can be realised
only by initialising the state of our system in some way that requires clairvoyant 
knowledge of the entire future of
its driving noise. 

In the framework presented in the previous section, such unphysical invariant measures
never arise, since our noise space does only contain information about the `past' of the noise.
Actually, given a skew-product of Markov processes as before, one can construct in a canonical way a random dynamical
system by taking $\Omega = \CW^\Z$, $\Theta$ the shift map, and $\bar\P$ the measure on
$\Omega = \CW^\Z$ corresponding to the law of the stationary Markov process with transition
probabilities $\CP$ and one-point distribution $\P$. The map $\bar \Phi$ is then given by
$\bar \Phi(x, \omega) = \Phi(x, \omega_0)$. With this correspondence, an invariant
measure for the skew-product yields an invariant measure for the corresponding
random dynamical system but, as the example given above shows, the converse is not
true in general.

\section{Ergodicity criteria for Markov semigroups}
\label{sec:ergodicity}

Consider a Markov transition kernel $\CP$ on some Polish space $\CX$.
Recall that an invariant measure $\mu$ for $\CP$ is said to be ergodic if the law of the
corresponding stationary process is ergodic for the shift map. It is a well-known fact that
if a Markov transition kernel has more than one invariant measure, then it must have at least two
of them that are mutually singular. Therefore, the usual strategy for
proving the uniqueness of the invariant measure for $\CP$ is to assume
that $\CP$ has two mutually singular invariant measures $\mu$ and $\nu$ and to
arrive at a contradiction.

This section is devoted to the presentation of some ergodicity criteria for
a Markov process on a general state space $\CX$ and to their extension to the
framework presented in Section~\ref{sec:framework}. If $\CX$ happens to be countable
(or finite), the transition probabilities are given by a transition matrix $P = (P_{ij})$
($P_{ij}$ being the probability of going from $i$ to $j$ in time $1$)
and there is a very simple characterisation of those transition probabilities that
can lead to at most one invariant probability measure. In a nutshell, ergodicity is implied
by the existence of one point which cannot be avoided by the dynamic:

\begin{proposition}
Let $P$ be a transition matrix. If there exists a state $j$ such that, for every $i$, $\sum_{n \ge 1} (P^n)_{ij} > 0$, then $P$ can have at most one invariant probability measure. Conversely, if $P$ has exactly one invariant probability measure, then there exists a state $j$ with the above property.
\end{proposition}

There is no such clean criterion available in the case of general state space,
but the following comes relatively close.
Recall that a Markov transition
operator $P$ over a Polish space $\CX$ is said to be 
\textit{strong Feller} if it maps the space of bounded measurable functions into the
space of bounded continuous functions. This is equivalent to the continuity
of transition probabilities in the topology of strong convergence of measures.
With this definition, one has the following criterion, of which a proof can be found
for example in \cite{DPZ96}:

\begin{proposition}\label{prop:disj}
Let $P$ be a strong Feller Markov transition operator on a Polish space $\CX$.
If $\mu$ and $\nu$ are two invariant measures for $P$ that are mutually singular,
then $\mathop{\mathrm{supp}} \mu \cap \mathop{\mathrm{supp}} \nu = \emptyset$.
\end{proposition}

This is usually used together with some controllability argument 
in the following way:

\begin{corollary}
Let $P$ be strong Feller.
If there exists $x$ such that $x$ belongs to the support of every invariant measure
for $P$, then $P$ can have at most one invariant measure.
\end{corollary}

\begin{remark}
The importance of the strong Feller property is that it allows to replace a measure-theoretical
statement ($\mu$ and $\nu$ are mutually singular) by a stronger topological statement (the topological
supports of $\mu$ and $\nu$ are disjoint) which is then easier to invalidate by a controllability
argument. If one further uses the fact that $\mu$ and $\nu$ are invariant measures, one can
actually replace the strong Feller property by the weaker asymptotic strong Feller property
\cite{HM06}, but we will not consider this generalisation here.
\end{remark}

A version with slightly stronger assumptions that however leads to a substantially stronger conclusion
is usually attributed to Doob and Khasminsk'ii:
\begin{theorem}\label{prop:DoobKha}
Let $P$ be a strong Feller Markov transition operator on a Polish space $\CX$.
If there exists $n \ge 1$ such that, for every open set $A \subset \CX$ and every $x \in \CX$,
one has $P^n(x,A) > 0$, then the measures $P^{m}(x, \cdot\,)$ and $P^{m}(y,\cdot\,)$ are
equivalent for every pair $(x,y) \in \CX^2$ and for every $m > n$. In particular, 
$P$ can have at most one invariant 
probability measure and, if it exists, it is equivalent to $P^{n+1}(x, \cdot\,)$ for every $x$.
\end{theorem}

These criteria suggest that we should look for a version of the strong Feller property
that is suitable for our context. Requiring $\CQ$ to be strong Feller is a very strong 
requirement which will not be fulfilled in many cases of interest. On the other hand,
since there is nothing like a semigroup on $\CX$, it is not clear \textit{a priori} how
the strong Feller property should be translated to our framework. 
On the other hand, the \textit{ultra Feller} property, that is the continuity of the transition
probabilities in the total variation topology is easier to generalise to our setting. 
Even though this property seems at first sight to be stronger than the strong Feller
property (the topology on probability measures induced by the total variation distance is strictly
stronger than the one induced by strong convergence), it turns out that the two
are `almost' equivalent. More precisely, if two Markov transition operators $P$ and $Q$ are strong Feller,
then $PQ$ is ultra Feller. Since this fact is not easy to find in the literature,
we will give a self-contained proof in Appendix~\ref{sec:ultrastrong} below.

One possible generalisation, and this is the one that we will retain here, is given by the following:
\begin{definition}\label{def:strongFeller}
A skew-product $(\CW, \P, \CP, \CX, \Phi)$ is said to be \textit{strong Feller} if there
exists a measurable map $\ell \colon \CW \times \CX^2 \to [0,1]$ such that, for $\P$-almost every $w$ one has $\ell(w,x,x) = 0$ for every $x$, and such that
\begin{equ}[e:defell]
\|\CS(x,w) - \CS(y,w)\|_\TV \le \ell(w,x,y)\;,
\end{equ}
for every $w\in \CW$ and every $x,y\in \CX$.

If we are furthermore in the setting of Section~\ref{sec:prod}, we assume that, 
for $\P$-almost every $w$, the map
$(w',x,y) \mapsto \ell(w\sqcup w',x,y)$ with $w' \in \CW_0$
is jointly continuous.

If we are not in that setting, we impose the stronger condition
that $\ell$ is jointly continuous. 
\end{definition}
A natural generalisation of the topological irreducibility used in Theorem~\ref{prop:DoobKha}
is given by
\begin{definition}\label{def:irreducible}
A skew-product $(\CW, \P, \CP, \CX, \Phi)$ is said to be \textit{topologically irreducible}
 if there exists $n \ge 1$ such that
\begin{equ}
\CQ^n(x,w;A \times \CW) > 0\;,
\end{equ}
for every $x \in \CX$, $\P$-almost every $w \in \CW$, and every open set $A \subset \CX$.
\end{definition}
According to these definitions, the example given in Section~\ref{sec:example} is
both strong Feller  and topologically irreducible. If $p \neq {1\over 2}$, it does however
have two distinct (even up to the equivalence relation $\simeq$) invariant measures. 
The problem is that in the non-Markovian case it is of course
perfectly possible to have two distinct ergodic invariant
measures for $\CQ$ that are such that their projections on $\CX$ are \textit{not} mutually singular.
This
shows that if we are aiming for an extension of a statement along the lines
of Theorem~\ref{prop:DoobKha}, we should impose some additional condition, which ideally should
always be satisfied for Markovian systems.

\subsection{Off-white noise systems}
\label{sec:offwhite}

In order to proceed, we consider the measure $\hat \P$ on $\CW^\Z$, which is the law of the
stationary process with transition probabilities $\CP$ and fixed-time law $\P$. We define
the coordinate maps $\Pi_i \colon \CW^\Z \to \CW$ in the natural way and the shift map $\Theta$
satisfying $\Pi_i \Theta w = \Pi_{i+1} w$.
We also define two natural $\sigma$-fields on $\CW^\Z$. The past, $\cP$, is defined as the
$\sigma$-field generated by the coordinate maps $\Pi_i$ for $i \le 0$. The future, $\cF$, is defined as the
$\sigma$-field generated by all the maps of the form $w \mapsto \Phi(x, \Pi_i w)$ for $x \in \CX$ and $i > 0$.
With these definitions, we see that the process corresponding to our skew-product is Markov (which 
is sometimes expressed by saying that the system is a `white noise system') if
$\cP$ and $\cF$ are independent under $\hat \P$.

A natural weakening of the Markov property is therefore given by:
\begin{definition}
A skew-product $(\CW, \P, \CP, \CX, \Phi)$ is said to be an \textit{off-white noise system}
if there exists a probability measure $\hat \P_0$ that is equivalent to $\hat \P$ and such that $\cP$ and
$\cF$ are independent under  $\hat \P_0$.
\end{definition}
\begin{remark}
The terminology ``off-white noise system'' is used by analogy on the one hand with ``white noise systems''
in the theory of random dynamical systems and on the other hand with ``off-white noise'' (or ``slightly coloured noise'')
as studied by Tsirelson in \cite{Tsi00Coloured,Tsi02White}.
\end{remark}

An off-white noise system behaves, as far as ergodic properties are concerned, pretty much like
a white noise (Markovian) system. This is the content of the following proposition:

\begin{theorem}\label{theo:offwhite}
Let $(\CW, \P, \CP, \CX, \Phi)$ be an off-white noise system and let $\mu$ and $\nu$ be two stationary 
measures for $\CQ$ such that $\CS \mu \perp \CS \nu$. Then their projections $\Pi_\CX^* \mu$
and $\Pi_\CX^* \nu$ onto the state space $\CX$ are also mutually singular.
\end{theorem}

\begin{proof}
Denote by $\hat \Phi\colon \CX \times \CW^\Z \to \CX^\N$ the solution map defined recursively by
\begin{equ}[e:defPhihat]
\bigl(\hat\Phi(x,w)\bigr)_0 = x\;,\qquad \bigl(\hat\Phi(x,w)\bigr)_n = \Phi\bigl(\bigl(\hat \Phi(x,w)\bigr)_{n-1}, \Pi_n w\bigr)\;.
\end{equ}
It follows from the construction that $\hat \Phi$ is $\CB(\CX) \otimes \cF$-measurable, where $\CB(\CX)$ denotes
the Borel $\sigma$-algebra of $\CX$. Denote now by $\mu_w$ and $\nu_w$ the disintegrations of $\mu$ and $\nu$
over $\CW$, that is the only (up to $\P$-negligible sets) measurable functions from $\CW$ to $\cM_1(\CX)$ such that
\begin{equ}
\mu(A \times B) = \int_B \mu_w(A) \P(dw)\;, \quad A\in \CB(\CX)\;,\quad B \in \CB(\CW)\;,
\end{equ}
and similarly for $\nu$. Using this, we construct measures $\hat \mu$ and $\hat \nu$ on
$\CX \times \CW^\Z$ by
\begin{equ}[e:defmuhat]
\hat \mu(A \times B) = \int_B \mu_{\Pi_0 w}(A) \hat \P(dw)\;, \quad A\in \CB(\CX)\;,\quad B \in \CB(\CW^\Z)\;.
\end{equ}
With these constructions, one has $\CS \mu = \hat \Phi^* \hat \mu$ and similarly for $\nu$. We also 
define $\hat \mu_0$ by the same expression as \eref{e:defmuhat} with $\hat \P$ replaced by $\hat \P_0$.
Since $\hat \P_0 \approx \hat \P$, one has $\hat \mu \approx \hat \mu_0$. However, since $\Pi_0$ is $\cP$-measurable
and since $\cP$ and $\cF$ are independent under $\hat \P_0$, one has $\hat \mu_0 \approx \Pi_\CX^*  \hat \mu_0
\otimes \hat \P_0$ when restricted to the $\sigma$-algebra $\CB(\CX) \otimes \cF$. This implies in particular that
\begin{equ}
\CS \mu = \hat \Phi^* \hat \mu \approx \hat \Phi^* \hat \mu_0 =   \hat \Phi^*(\Pi_\CX^*  \hat \mu_0 \otimes \hat \P_0)
\approx  \hat \Phi^*(\Pi_\CX^*  \mu \otimes \hat \P)\;,
\end{equ}
which concludes the proof.
\end{proof}

\begin{remark}
Our definition of off-white noise is slightly more restrictive than the one in \cite{Tsi02White}. Translated to our present setting,
Tsirelson defined $\cF_n$ as the $\sigma$-field generated by all the maps of the 
form $w \mapsto \Phi(x, \Pi_i w)$ for $x \in \CX$ and $i > n$ and a noise was called ``off-white'' if there
exists $n \ge 0$ and $\hat \P_0 \approx \hat \P$ such that $\cP$ and $\cF_n$ are independent under $\hat \P_0$.

With this definition, one could expect to be able to obtain a statement similar to Theorem~\ref{theo:offwhite}
with the projection on $\CX$ of $\mu$ (and $\nu$) replaced by the projection on the first $n+1$ copies
of $\CX$ of the solution $\CS \mu$. Such a statement is {\em wrong}, as can be seen again by the example
from Section~\ref{sec:example}. It is however true if one defines $\cF_n$ as the (larger) $\sigma$-algebra generated
by all maps of the form $w \mapsto \bigl(\hat \Phi(x,w)\bigr)_i$ for $x \in \CX$ and $i > n$.
\end{remark}

A consequence of this theorem is the following equivalent of Proposition~\ref{prop:disj}:

\begin{proposition}
Let $(\CW, \P, \CP, \CX, \Phi)$ be an off-white noise system which is strong Feller in the sense
of Definition~\ref{def:strongFeller}. 
If there exists $x \in \CX$ such that
$x \in \mathop{\mathrm{supp}} \Pi_\CX^* \mu$ for every stationary measure $\mu$ of $\CQ$, then
there can be at most one such measure, up to the equivalence relation \eref{e:equivalence}.
\end{proposition}

\begin{proof}
Assume by contradiction that there 
exist two distinct invariant measures $\mu$ and $\nu$.
For simplicity, denote $\mu_\CX = \Pi_\CX^* \mu$ and similarly for $\nu$ and let $x$ be an element
from the intersection of their supports (such an $x$ exists by assumption). 
We can assume furthermore without any loss of generality
that $\CS \mu \perp \CS \nu$.

Define, with the same notations as in the proof of Theorem~\ref{theo:offwhite}, 
\begin{equ}
\hat \CS(x;\cdot\,) =  \hat \Phi^* \bigl(\delta_x \otimes \hat \P\bigr) = \int_{\CW} \CS(x,w;\cdot\,)\,\P(dw)
\end{equ}
and note that one has, as before,
$\CS \mu \approx \int \hat \CS(x;\cdot\,)\,\mu_\CX(dx)$ and similarly for $\nu$. Since $\CS \mu \perp \CS\nu$,
this shows that there exists a measurable set $A$ such that $\CS(x,w;A) = 0$ for $\mu_\CX \otimes \P$-almost
every $(x,w)$ and $\CS(x,w;A) = 1$ for $\nu_\CX \otimes \P$-almost every $(x,w)$. Define a function $\delta \colon \CW \to \R_+$ by
\begin{equ}
\delta(w) = \inf \{\delta \,:\, \exists x_0 \;\hbox{with $d(x_0,x) < \delta$ and $\CS(x_0,w;A) = 0$}\}\;,
\end{equ}
were $d$ is any metric generating the topology of $\CX$. Since $x$ belongs to the support of $\mu_\CX$ one
must have $\delta(w) = 0$ for $\P$-almost every $w$. Since we assumed that $y \mapsto \CS(y,w;\cdot\,)$
is continuous, this implies that $\CS(x,w;A) = 0$ for $\P$-almost every $w$. Reversing the roles of $\mu$ and $\nu$, one
arrives at the fact that one also has $\CS(x,w;A) = 1$ for $\P$-almost every $w$, which is the contradiction we were looking for.
\end{proof}

\begin{remark}
It follows from the proof that it is sufficient to assume that the map $x \mapsto \CS(x,w)$ is continuous
in the total variation norm for $\P$-almost every $w$. 
\end{remark}

\subsection{Another quasi-Markov property}

While the result in the previous section is satisfactory in the sense that it shows a nice correspondence between 
results for Markov processes and results for off-white noise systems, it covers only a very restrictive class of systems.
For example, in the case of continuous time,
neither fractional noise (the derivative of fractional Brownian motion) nor the Ornstein-Uhlenbeck
process fall into this class. It is therefore natural to look for weaker conditions that still allow to obtain statements
similar to Theorem~\ref{prop:DoobKha}. The key idea at this stage is to make use of the topology of $\CW$
which has not been used in the previous section. This is also the main conceptual difference between
the approach outlined in this article and the approach used by the theory of random dynamical
systems.

In the previous section, we made use of the fact that for off-white noise systems, one has
$\CS(x,w;\cdot\,) \approx \CS(x,w';\cdot\,)$ for every $x$ and every pair $(w,w')$ in a set of full $\P$-measure.
We now consider the set
\begin{equ}[e:defLambda]
\Lambda = \{(w,w') \in \CW^2\,:\,\CS(x,w;\cdot\,) \approx \CS(x,w';\cdot\,)\}\;,
\end{equ}
and we require that the dynamic on $\CW$ is such that one can construct couplings that hit $\Lambda$
with positive probability. If we think of the driving noise to be some Gaussian process, the set $\Lambda$ typically
consists of pairs $(w,w')$ such that the difference $w - w'$ is sufficiently `smooth'.

Recall that a \textit{coupling} between two probability measures $\mu$ and $\nu$ is a measure $\pi$
on the product space such that its projections on the two factors are equal to $\mu$ and $\nu$ respectively
(the typical example is $\pi = \mu \otimes \nu$ but there exist in general many different couplings
for the same pair of measures). Given two positive measures $\mu$ and $\nu$, we say that $\pi$
is a \textit{subcoupling} for $\mu$ and $\nu$ if the projections on the two factors are \textit{smaller} than
$\mu$ and $\nu$ respectively. With this definition at hand, we say that:

\begin{definition}\label{def:quasiMarkov}
A skew-product $(\CW, \P, \CP, \CX, \Phi)$ is said to be \textit{quasi-Mar\-ko\-vian}
if, for any two open sets $U,V \subset \CW$ such that $\min\{\P(U),\P(V)\} > 0$, there exists a measurable map $w \mapsto \CP^{U,V}(w,\cdot\,)
\in \cM_+(\CW^2)$ such that:
\begin{itemize}
\item[i)] For $\P$-almost every $w$, the measure $\CP^{U,V}(w,\cdot\,)$ is a subcoupling for $\CP(w,\cdot\,)|_U$ and $\CP(w,\cdot\,)|_V$.
\item[ii)] Given $\Lambda$ as in \eref{e:defLambda}, one has $\CP^{U,V}(w,\Lambda) = \CP^{U,V}(w,\CW^2)$ for $\P$-almost every $w$.
\end{itemize}
\end{definition}

\begin{remark}
If we are in the setting of Section~\ref{sec:prod}, this is equivalent to considering for $U$ and $V$ open sets in $\CW_0$
and replacing every occurrence of $\CP$ by $\bar \CP$. The set $\Lambda$ should then be replaced by the set
\begin{equ}
\bar \Lambda = \{(w_0,w_0') \in \CW_0^2 \,:\, \CS(x, w \sqcup w_0;\cdot\,) \approx \CS(x, w \sqcup w_0';\cdot\,) \; \hbox{for $\P$-almost every $w$}\}
\end{equ}
\end{remark}

\begin{remark}\label{rem:negligible}
In general, the transition probabilities $\CP$ only need to be defined up to a $\P$-negligible set.
In this case, the set $\Lambda$ is defined up to a set which is negligible with respect to any 
coupling of $\P$ with itself. In particular, this shows that the ``quasi-Markov'' property from
Definition~\ref{def:quasiMarkov} does not depend on the particular choice of $\CP$.
\end{remark}

With these definitions, we have the following result:

\begin{theorem}
Let $(\CW, \P, \CP, \CX, \Phi)$ be a quasi-Markovian skew-product which is strong Feller in the sense
of Definition~\ref{def:strongFeller} and topologically irreducible in the sense of Definition~\ref{def:irreducible}.
Then, it can have at most one invariant measure, up to the equivalence relation \eref{e:equivalence}.
\end{theorem}

\begin{proof}
Under slightly more restrictive assumptions, this is the content of \cite[Theorem~3.10]{HaiOha07}.
It is a tedious but
rather straightforward task to go through the proof and to check that the arguments still hold under
the weaker assumptions stated here.
\end{proof}

\subsection{Discussion}

The insight that we would like to convey with the way of exposing the previous two subsections is the following.
If one wishes to obtain a statement of the form ``strong Feller + irreducible + quasi-Markov $\Rightarrow$ uniqueness of the stationary measure'', one should balance the regularity of $\ell$, defined in \eref{e:defell}, as a function
of $w$, with the
class of sets $U$ and $V$ used in Definition~\ref{def:quasiMarkov}. This in turn is closely related to the
size of the set $\Lambda$ from \eref{e:defLambda}. The larger $\Lambda$ is, the larger the admissible class of sets
in Definition~\ref{def:quasiMarkov}, and the lower the regularity requirements on $\ell$.

The off-white noise case corresponds to the situation where $\Lambda = \CW^2$. This in turn shows that one could take
for $U$ and $V$ any two measurable sets and $\CP^{U,V}(w,\cdot\,) = \CP(w,\cdot\,)|_U \otimes \CP(w,\cdot\,)|_V$. 
Accordingly, there is no regularity requirement (in $w$) on $\ell$,  except for it being measurable.

In the case of Section~\ref{sec:prod}, the transition probabilities $\CP$ have a special structure in the sense that
$\CP(w, A) = 1$ for $A = w \sqcup \CW_0$. This implies that one can take for $U$ and $V$ any measurable set 
that is such that, if we decompose $\CW$ according to $\CW \approx \CW \times \CW_0$, the ``slices'' of $U$ and $V$
in $\CW_0$ are $\P$-almost surely open sets. The corresponding regularity requirement on $\ell$ is that
the map $(x,y,w') \mapsto \ell(w \sqcup w', x, y)$ is jointly continuous for $\P$-almost every $w$.

Finally, if we do not assume any special structure on $\Lambda$ or $\CP$, we take for $U$ and $V$ arbitrary open sets in $\CW$.
In this case, the corresponding regularity requirement on $\ell$ is that it is jointly continuous in all of its arguments.

\section{The Gaussian case}
\label{sec:Gauss}

In this section, we study the important particular case of Gaussian noise. We place ourselves
in the framework of Section~\ref{sec:prod} and we choose $\CW_0 = \R$, so that $\CW = \R^{\Z_-}$. We furthermore
assume that the measure $\P$ is centred, stationary, and Gaussian with covariance
$C$ and spectral measure $\mu$. In other words, we define $\mu$ as the (unique) finite
Borel measure on $[-\pi,\pi]$ such that
\begin{equ}[e:defspectral]
C_n = \int_\CW w_k w_{k+n} \, \P(dw) = {1\over 2\pi} \int_{-\pi}^\pi e^{inx}\mu(dx)\;,
\end{equ}
holds for every $n \ge 0$.
A well-known result by Maruyama, see \cite{Mar49} or the textbook
 \cite[Section~3.9]{DM76Gauss}, states that $\P$ is ergodic for the shift map if and only if the measure
$\mu$ has no atoms. As in Section~\ref{sec:prod}, denote by $\bar \CP\colon \CW \to \cM_1(\R)$ the
corresponding regular conditional probabilities.

Since regular conditional probabilities of Gaussian measures are again Gaussian \cite{Bog98Gauss}, one has
\begin{lemma}\label{lem:cond}
There exists $\sigma \ge 0$ and a $\P$-measurable linear functional $\mean \colon \CW \to \R$ 
such that, $\P$-almost surely, the measure  $\bar \CP(w,\cdot\,)$ is Gaussian with mean $\mean(w)$ and
variance $\sigma^2$.
\end{lemma}
This however does not rule out the case where $\sigma = 0$.
The answer to the question of when $\sigma \neq 0$ is given
by the following classical result in linear prediction theory \cite{Szego20,HS60}:
\begin{theorem}
Decompose $\mu$ as $\mu(dx) = f(x)\, dx + \mu_s(dx)$ with $\mu_s$ singular with respect to 
Lebesgue's measure. Then, one has
\begin{equ}[e:exprsigma]
\sigma^2 = \exp \Bigl({1\over 2\pi}\int \log f(x)\, dx\Bigr)\;,
\end{equ}
if the expression on the right hand side makes sense and $\sigma = 0$ otherwise.
\end{theorem}
If $\sigma = 0$, all the randomness is contained in the remote past of the noise and no
new randomness comes in as time evolves. We will therefore always assume that 
$\mu$ is non-atomic and that 
$\sigma^2 > 0$. Since in that case all elements of $\CW$ with only finitely many non-zero
entries belong to the reproducing kernel of $\P$ (see Section~\ref{sec:Gaussian} below for the
definition of the reproducing kernel of a Gaussian measure
and for the notations that follow), the linear functional $\mean$ can be chosen
such that, for every $n>0$, $\mean(w)$ is jointly continuous in $(w_{-n},\ldots, w_0)$ for $\P$-almost
every $(\ldots, w_{-n-2}, w_{-n-1})$, see \cite[Sec.~2.10]{Bog98Gauss}.

We will denote by $\hat \P$ the Gaussian measure on $\hat \CW = \R^\Z$ with correlations given by
the $C_n$. We denote its covariance operator again by $C$.
The measure $\hat \P$ is really the same as the measure $\hat \P$ defined in 
Section~\ref{sec:offwhite} if we make the necessary identification of $\hat \CW$ with a subset of $\CW^\N$, this
is why we use the same notation without risking confusion. 
We also introduce the equivalents to the two $\sigma$-algebras $\cP$ and $\cF$.
We interpret them as $\sigma$-algebras on $\R^\Z$, so that $\cP$ is the $\sigma$-algebra generated by the $\Pi_n$
with $n \le 0$ and $\cF$ is generated by the $\Pi_n$ with $n > 0$.
(Actually, $\cF$ could be slightly smaller than that in general, but we do not want to restrict ourselves to one 
particular skew-product, and so we simply take for $\cF$ the smallest choice which contains all `futures' for all
possible choices of $\Phi$ as in Section~\ref{sec:prod}.)

It is natural to split $\hat \CW$ as $\hat \CW = \CW_- \oplus \CW_+$ where $\CW_- \approx \CW$ is the span of the
images of the $\Pi_n$ with $n \le 0$ and similarly for $\CW_+$.
We denote by $\CH$ the reproducing kernel Hilbert space of $\hat \P$.
Recall that via the map $\hat \CW^* \ni \Pi_n \mapsto e^{inx}$ and the inclusion $\hat \CW^* \subset \CH$, 
one has the isomorphism $\CH \approx L^2(\mu)$ with
$\mu$ as in \eref{e:defspectral}, see for example \cite{DM76Gauss}. 
Following the construction of Section~\ref{sec:Gaussian}, we
see that $\hat \CH_-^p$ is given by the closure in $\CH$ of the span of $e^{inx}$ for $n \le 0$ and similarly for $\hat \CH_-^p$.
Denote by $P_\pm$ the orthogonal projection from $\hat \CH_-^p$ to $\hat \CH_+^p$ and by $P$ the
corresponding operator from $\CH_-^p$ to $\CH_+^p$.

With all these preliminaries in place, an immediate consequence of Proposition~\ref{prop:conditional} is:
\begin{proposition}\label{prop:LambdaGauss}
Let $\hat P \colon \CW_- \to \CW_+$ be the $\P$-measurable extension of $P$. Then, the set
$\Lambda$ is equal (up to a negligible set in the sense of Remark~\ref{rem:negligible}) to 
$\{(w,w') \,:\, \hat P (w - w') \in \CH_+^c\}$. In particular, it always contains the set 
$\{(w,w') \,:\, w - w' \in \hat \CH \cap \CW_-\}$.
\end{proposition}

\begin{proof}
The first statement follows from the fact that, by \eref{e:dis}, $\CH_+^c$ is the reproducing kernel space of 
the conditional probability of $\hat \CP$, given the past $\cP$, and $P(w-w')$ yields the shift between the
conditional probability given $w$ and the conditional probability given $w'$.
The second statement follows from the fact that $P$ extends to a bounded
operator from $\CH_-^c$ to $\CH_+^c$.
\end{proof}

\subsection{The quasi-Markov property}

We assume as above that $\CW_0 = \R$ and that $\P$ is a stationary Gaussian measure
with spectral measure $\mu$. We also write as before $\mu(dx) = f(x)\,dx + \mu_s(dx)$.
The main result of this section is that the quasi-Markov property introduced in Section~\ref{sec:ergodicity}
can easily be read off from the behaviour of the spectral measure $\mu$:
\begin{theorem}\label{thm:QM}
A generic random dynamical system as above is quasi-Markovian if and only if 
$f$ is almost everywhere positive and $\int_{-\pi}^\pi {1 \over f(x)}\,dx$ is finite.
\end{theorem}

\begin{proof}
Let $e_n$ be the `unit vectors' defined by $\Pi_m e_n = \delta_{mn}$. Then the condition of
$\int_{-\pi}^\pi {1 \over f(x)}\,dx$ being finite is equivalent to $e_n$ belonging to the reproducing kernel
of $\hat \P$, a classical result dating back to Kolmogorov, see also \cite[p.~83]{GreRos57}.

To show that the condition is sufficient, denote by $D_w(x)$ the (Gaussian) density of $\bar \CP(w,\cdot\,)$ with respect to 
Lebesgue measure on $\R$. Given any two open sets $U$ and $V$ in $\R$, we can find some $x$, $y$, and $r> 0$ such that
$\CB(x,r) \subset U$ and $\CB(y,r) \subset V$. Take then for $\bar \CP^{U,V}(w,\cdot\,)$ 
the push-forward under the map
$z \mapsto (z, z+y-x)$ of the measure with density
$z \mapsto \min\{D_w(z), D_w(z + y - x)\}$ with respect to Lebesgue measure. Since, by Proposition~\ref{prop:LambdaGauss}, $\Lambda$ contains all pairs $(w,w')$ which
differ by an element of $\hat \CH \cap \CW_-$ and since the condition of the theorem is precisely what is
required for $e_0$ to belong to $\hat \CH$, this shows the sufficiency of the condition.

To show that the condition is necessary as well, suppose that it does not hold and take for example 
$U = (-\infty,0)$ and $V = [0,\infty)$. Since we have the standing assumption that $\sigma^2 > 0$ 
with $\sigma^2$ as in \eref{e:exprsigma},
one has $\bar \CP(w,U) > 0$ for $\P$-almost every $w$ and similarly for $V$. Assume by contradiction
that the system is quasi-Markovian, so we can construct a measure on $\R^\Z \times \R^\Z$ in the following way.
Define $\CW_+$ as before and define $\CW_-$ as the span of $\Pi_n$ for $n < 0$ so that $\CW = \CW_- \oplus \CW_0 \oplus \CW_+$.
Let $P_\pm \colon \CW_- \oplus \CW_0 \to \cM_1(\CW_+)$ be the conditional probability of $\hat \P$ given 
$\CW_- \oplus \CW_0$. 
Let $\bar \CP^{U,V}\colon \CW_- \to \cM_1(\CW_0)$ be as in Definition~\ref{def:quasiMarkov} and construct a 
measure $M$ on $ \CW_-^2 \times  \CW_+^2$  by
\begin{equs}
M(A_1 \times A_2 \times B_1 \times B_2) &= \int_{A_1 \cap A_2} \int_{\bar \Lambda} P_\pm(w_- \sqcup w_0, B_1) P_\pm(w_- \sqcup w_0', B_2)  \\
&\quad \times \bar \CP^{U,V}(w_-, dw_0 dw_0') \, \P(dw_-)\;.
\end{equs}
This measure has the following properties:
\begin{itemize}
\item[1.] By the properties of $\bar \CP^{U,V}$ and by the definition of $P_\pm$, it
 is a subcoupling for the projection of $\hat \P$ on $\CW_- \times \CW_+$ with itself, and it
 is not the trivial measure.
\item[2.] Denote by $M_1$ and $M_2$ the projections $M_1$ and $M_2$ on the two copies of $\CW_1 \times \CW_+$.
Since $P_\pm(w_- \sqcup w_0, \cdot) \approx P_\pm(w_- \sqcup w_0', \cdot)$ for $\P$-almost every $w$ 
and for every pair $(w_0, w_)') \in \bar \Lambda$, one has $M_1\approx M_2$.
\end{itemize}

On the other hand, since $e_0$ does not belong to the reproducing kernel
of $\hat \P$ by assumption, there exists a $\hat \P$-measurable linear map $\mean\colon \CW_- \times \CW_+ \to \CW_0$ such
that the identity $w_0 = \mean(w_-, w_+)$ holds for $\hat \P$-almost every triple $(w_-, w_0, w_+)$, see
Proposition~\ref{prop:conditional}. Denote by $A$ the preimage of $U$ under $\mean$ in $\CW_- \times \CW_+$
and by $A^c$ its complement. Then one has $M_1(A^c) = M_2(A) = 0$, which contradicts property 2 above.
\end{proof}

Note that although the condition of this theorem is easy to read off from the spectral measure,
it is in general not so straightforward to read off from the behaviour of the
correlation function $C$. In particular, it does \textit{not} translate into a decay condition
of the coefficients $C_n$. Take for example the case
\begin{equ}
C_n = \left\{\begin{array}{rl} 2 & \text{if $n = 0$,} \\ 1 & \text{if $n = 1$,} \\ 0 & \text{otherwise.} \end{array}\right.
\end{equ}
This can be realised for example by taking for $\xi_n$ a sequence of i.i.d.\ normal Gaussian random variables
and setting $W_n = \xi_n + \xi_{n+1}$. We can check that one has, for every $N > 0$, the identity
\begin{equ}
W_0 = {1\over N} \sum_{n=1}^N (-1)^n \bigl(\xi_{n+1} + \xi_{-n}\bigr) - {1\over N}\sum_{n=1}^N (-1)^n(N + 1-n)\bigl(W_n + W_{-n}\bigr) \;.
\end{equ}
Since the first term converges to $0$ almost surely by the law of large numbers, it follows that
one has the almost sure identity
\begin{equ}
W_0 = - \lim_{N\to \infty} {1\over N}\sum_{n=1}^N (-1)^n(N + 1-n)\bigl(W_n + W_{-n}\bigr) \;,
\end{equ}
which shows that $W_0$ can be determined from the knowledge of the $W_n$ for $n \neq 0$.
In terms of the spectral measure, this can be seen from the fact that $f(x) = 1 + \cos(x)$, so that $1/f$ has
a non-integrable singularity at $x=\pi$. This also demonstrates that there are cases in which the reproducing kernel
of $\P$ contains all elements with finitely non-zero entries, even though the reproducing kernel
of $\hat \P$ contains no such elements.

\subsection{The strong Feller property}

It turns out that in the case of discrete stationary Gaussian noise, the quasi-Markov and the strong Feller
properties are very closely related. In this section, we assume that we are again in the framework of Section~\ref{sec:prod},
but we take $\CW_0 = \R^d$ and we assume that the driving noise consists of $d$ independent stationary Gaussian
sequences with spectral measures satisfying the condition of Theorem~\ref{thm:QM}.

We are going to derive a criterion for the strong Feller property for the Markovian case where the driving noise consists of $d$ independent  sequences of 
i.i.d.\ Gaussian random variables and we will see that this criterion still works 
in the quasi-Markovian case.

It will be convenient for the purpose of this section
to introduce the Fr\'echet space $L^\Gamma(\R^d)$ consisting of measurable functions
$f \colon \R^d \to \R$ such that the norms $\|f\|_{\gamma,p}^p = \int f^p(w) e^{-\gamma |w|^2} \,dw$ are finite
for all $\gamma > 0$ and all $p \ge 1$. For example, since these norms are increasing in $p$ and and decreasing in $\gamma$, 
$L^\Gamma(\R^d)$ can be endowed with the distance
\begin{equ}
d(f,g) = \sum_{p=1}^\infty \sum_{n=1}^\infty 2^{-n-p}\bigl(1 \wedge \|f-g\|_{{1\over n}, p}\bigr)\;.
\end{equ}

With this notation, we will say that a function $g\colon \R^n \times \R^d \to \R^m$ belongs to $\CC^{0,\Gamma}(\R^n \times \R^d)$ if, for every $i \in \{1,\ldots,m\}$, the map $x \mapsto g_i(x,\cdot\,)$ is continuous from $\R^n$ to $L^\Gamma(\R^d)$.

Given a function $\Phi\colon \R^n \times \R^d \to \R^n$ with elements of $\R^n$ denoted by $x$ and elements of $\R^d$ denoted by $w$, we also define the ``Malliavin covariance matrix'' of $\Phi$ by
\begin{equ}
M^\Phi_{ij}(x,w) = \sum_{k=1}^d \d_{w^k}\Phi_i(x,w)\, \d_{w^k}\Phi_j(x,w)\;.
\end{equ}
With these notations, we have the following criterion:

\begin{proposition}\label{prop:SF}
Let $\Phi \in \CC^2(\R^n \times \R^d ; \R^n)$ be such that the derivatives $D_w \Phi$, $D_w D_x \Phi$, and $D_w^2 \Phi$ all belong to $\CC^{0,\Gamma}(\R^n \times \R^d)$. Assume furthermore that 
$M^\Phi_{ij}$ is invertible for Lebesgue-almost every $(x,w)$ and that $(\det  M^\Phi_{ij})^{-1}$ belongs
to $\CC^{0,\Gamma}(\R^n \times \R^d)$. Then, the Markov semigroup over $\R^d$ defined by
\begin{equ}
\bigl(\CP f\bigr)(x) = \int_{\R^d} f(\Phi(x,w))\, \Gamma(dw)\;,
\end{equ} 
where $\Gamma$ is an arbitrary non-degenerate Gaussian measure on $\R^d$, has the strong Feller
property.
\end{proposition}

\begin{proof}
Take a function $f \in \CC_0^\infty(\R^n)$ and write (in this proof we use Einstein's convention of summation over repeated indices):
\begin{equ}[e:dP]
\bigl(\d_i \CP f \bigr)(x) = \int_{\R^d} \d_j f(\Phi(x,w))\d_{x_i} \Phi_j(x,w)\, \Gamma(dw)\;.
\end{equ}
At this point, we note that since we assumed $M^\Phi$ to be invertible, one has for every pair $(i,j)$ 
the identity
\begin{equ}[e:ibp]
\d_{x_i} \Phi_j (x,w) = \d_{w^m} \Phi_j (x,w)\Xi_{mi}(x,w) \;.
\end{equ}
where
\begin{equ}
\Xi_{mi} =  \d_{w^m} \Phi_k (x,w)  \bigl(M^\Phi(x,w)\bigr)^{-1}_{k \ell} \d_{x_i} \Phi_\ell (x,w)
\end{equ}
This allows to integrate \eref{e:dP} by parts, yielding
\begin{equ}
\bigl(\d_i \CP f \bigr)(x) =  -\int_{\R^d} f(\Phi(x,w)) \bigl(\d_{w^m} - (Q w)^m \bigr)\Xi_{mi}(x,w)\, \Gamma(dw)\;,
\end{equ}
where $Q$ is the inverse of the covariance matrix of $\Gamma$. Our assumptions then ensure the existence of a continuous function $K \colon \R^n \to \R$ such that $|\bigl(\d_i \CP f \bigr)(x)| \le K(x) \sup_y |f(y)|$ which, by
a standard approximation argument \cite[Chapter~7]{DPZ96}, is sufficient for the strong Feller property to hold.
\end{proof}

\begin{remark}
We could easily have replaced $\R^n$ by an $n$-dimensional Riemannian manifold with the obvious
changes in the definitions of the various objects involved.
\end{remark}

\begin{remark}
Just as in the case of the H\"ormander condition for the hypoellipticity of a second-order differential operator,
the conditions given here are not far from being necessary. Indeed, if $M^\Phi(x,\cdot\,)$ fails to be invertible on 
some open set in $\R^d$, then the image of this open set under $\Phi(x,\cdot\,)$ will be a set of dimension $n' < n$.
In other words, the process starting from $x$ will stay in some subset of lower dimension $n'$ 
with positive probability, so that
the transition probabilities will not have a density with respect to the Lebesgue measure.
\end{remark}

\begin{remark}
Actually, this condition gives  quite a bit more than the strong Feller property, since it gives local Lipschitz continuity of the 
transition probabilities in the total variation distance with local Lipschitz constant $K(x)$.
\end{remark}

We now show that if we construct a skew-product from $\Phi$ and take as driving noise $d$ independent copies
of a stationary Gaussian process with a covariance structure satisfying the assumption of Theorem~\ref{thm:QM},
then the assumptions of Proposition~\ref{prop:SF} are sufficient to guarantee that it also satisfies the strong Feller
property in the sense of Definition~\ref{def:strongFeller}. We have indeed that:

\begin{theorem}\label{thm:SFNM}
Let $\CW_0 = \R^d$, $\CW = \CW_0^{\Z_-}$, let $\Phi\colon \R^n \times \CW_0 \to \R^n$ satisfy the assumptions of Proposition~\ref{prop:SF}, 
and let $\P \in \cM_1(\CW)$ be a Gaussian measure such that there exist measures $\mu_1 ,\ldots, \mu_d$
with
\begin{equ}
\int_\CW w_n^j w_{n+k}^m \, \P(dw) = {\delta_{jm} \over 2\pi} \int_{-\pi}^{\pi} e^{ik x} \mu_j(dx)\;.
\end{equ}
Then, if the absolutely continuous part of each of the $\mu_j$ satisfies the condition of Theorem~\ref{thm:QM},
the skew-product $(\CW, \P, \CP, \R^n, \Phi)$ has the strong Feller property.
\end{theorem}

\begin{proof}
Let $\mean$ be as in Lemma~\ref{lem:cond}, let $x \in \R^n$, and let $w \in \CW$ such that $\mean(w) < \infty$.
We want to show that there exists a continuous function 
$\bar K$ depending continuously on $x$ and on $\mean(w)$ such that $\CS(x,w;\cdot\,)$ is locally
Lipschitz continuous in the total variation distance (as a function of $x$) with local 
Lipschitz constant $\bar K(x,\mean(w))$.
Since we assumed from the beginning that $\sigma^2 > 0$, where $\sigma$ is defined as in 
\eref{e:exprsigma}, we know that the set of elements in $\CW$ with only finitely many non-zero
coordinates belongs to the reproducing kernel space of $\P$. Since, by 
\eref{e:boundP}, the map $\mean$ is bounded from the reproducing kernel space of $\P$
into $\R$, $\mean(w)$ depends continuously on each of the coordinates of $w$ and so the
assumptions of Definition~\ref{def:strongFeller} are verified.

It remains to construct $\bar K$. This will be done in a way that is almost identical to the proof of
Proposition~\ref{prop:SF}. Take a bounded smooth test function $f \colon (\R^n)^\N \to \R$ 
which depends only on its first $N$ coordinates and consider
the function $\bigl(\CS f\bigr)(x,w)$ defined by
\begin{equ}
\bigl(\CS f\bigr)(x,w) = \int_{(\R^n)^\N} f(y) \CS(x,w;dy)\;.
\end{equ}
Consider now the splitting $\hat \CW = \CW_- \oplus \CW_+$, as well as the measurable
linear map $\hat P$ and the space $\CH_+^c$ introduced for the statement of
Proposition~\ref{prop:LambdaGauss} (note that $\hat P$ relates to $\mean$ via
$(\hat Pw)_0 = \mean(w)$). 
Denote furthermore by $\P_+$ the Gaussian measure on
$\CW_+$ with reproducing kernel space $\CH_+^c$. With these notations at hand, we have the
expression
\begin{equ}
\bigl(\CS f\bigr)(x,w) = \int_{\CW_+} f \bigl(\hat \Phi(x,\tilde w + \hat P w) \bigr) 
\P_+(d\tilde w)\;,
\end{equ}
where we denoted by $\hat \Phi \colon \CX \times \CW_+
\to \CX^\N$ the map defined in \eref{e:defPhihat}.
We see that, as in \eref{e:ibp}, one has the identity (again, summation over repeated indices is implied):
\begin{equ}
\d_{x_i} \hat \Phi(x,\tilde w) = \d_{\tilde w_0^m} \hat \Phi(x,\tilde w)\, \Xi_{mi}(x,\tilde w_0)\;,
\end{equ}
where the function $\Xi$ is exactly the same as in \eref{e:ibp}. At this point, since the `coordinate vectors' 
$e_0^m$ belong to the reproducing kernel of $\hat \P$, and therefore also of $\P_+$ by
\eref{e:carCHc}, we can integrate by parts against the Gaussian measure $\P_+$ 
\cite[Theorem~5.1.8]{Bog98Gauss} to obtain
\begin{equs}
\d_i \CS f(x,w) &= - \int_{\CW_+} f(\hat \Phi(x,\tilde w + \hat Pw)) \d_{\tilde w_0^m} \Xi_{mi}(x,\tilde w_0 + \mean(w))\, \P_+(d\tilde w) \\
&\quad + \int_{\CW_+} f(\hat \Phi(x,\tilde w+ \hat Pw)) \Xi_{mi}(x,\tilde w_0 + \mean(w)) e_0^m(\tilde w)\,\P_+(d\tilde w)\;.
\end{equs}
Here, we made an abuse of notation and interpreted $e_0^m$ as a measurable linear functional
on $\CW_+$, via the identification \eref{e:ident}. Since $\int |e_0^m(\tilde w)|^2 \P_+(d\tilde w) < \infty$ by assumption and since the law of $\tilde w_0$ under $\P_+$ is centred Gaussian with variance $\sigma^2$, this concludes the proof.
\end{proof}

\subsection{The off-white noise case}

The question of which stationary Gaussian sequences correspond to off-white noise
systems was solved by Ibragimov and Solev in the seventies, see \cite{IR78Gaussian} 
and also \cite{Tsi02White}. It turns out that the correct criterion is:

\begin{theorem}
The random dynamical system is off-white if and only if the spectral measure $\mu$ has a density $f$ with respect to Lebesgue measure and $f(\lambda) = \exp \phi(\lambda)$ for some function
$\phi$ belonging to the fractional Sobolev space $H^{1/2}$.
\end{theorem}

\begin{remark}
It follows from a well-known result by \cite{Tru67}, later extended in \cite{Str72}, that any
function $\phi \in H^{1/2}$ satisfies $\int \exp(\phi^2(x))\, dx < \infty$. In particular, this shows that
the condition of the previous theorem is therefore much stronger than the condition of Theorem~\ref{thm:QM}
which is required for the quasi-Markov property.
\end{remark}

As an example, the (Gaussian) stationary autoregressive process, which has covariance structure
$C_n = \alpha^n$ does have the quasi-Markov property since its spectral measure 
has a density of the form
\begin{equ}
\mu(dx) = {1-\alpha^2 \over 1+\alpha^2 - 2\alpha \cos(x)}\,dx\;,
\end{equ}
which is smooth and bounded away from the origin. However, if we take a sequence $\xi_n$ of
i.i.d.\ normal random variables and define a process $X_n$ by
\begin{equ}
X_n = \sum_{k = 1}^\infty k^{-\beta} \xi_{n-k}\;,
\end{equ}
for some $\beta > 1/2$, then $X_n$ does still have the quasi-Markov property, but it 
is not an off-white noise.

\section[Strong and ultra Feller property]{Appendix A: Equivalence of the strong and ultra Feller properties}
\label{sec:ultrastrong}

In this section, we show that, even though the ultra Feller property
seems at first sight to be stronger than the strong Feller property, the composition 
of two Markov transition kernels 
satisfying the strong Feller property always satisfies the ultra Feller property. 
This fact had already been pointed out in \cite{DM83} but had been overlooked by
a large part of the probability community until Seidler `rediscovered' it in 2001 \cite{Sei02}.
We take the opportunity to give an elementary proof of this fact.
Its structure is based on the notes by Seidler, but we take advantage
of the simplifying fact that we only work with Polish spaces.

We introduce the following definition:

\begin{definition}
A Markov transition kernel $P$ over a Polish space $\CX$ satisfies the \textit{ultra Feller}
property if the transition probabilities $P(x,\cdot\,)$ are continuous in the total variation
norm.
\end{definition}

Recall first the following well-known fact of real analysis, see for example 
\cite[Example~IV.9.3]{Yos95}:
\begin{proposition}
For any measure space $(\Omega,\CF,\lambda)$ such that $\CF$ is countably generated and any $p \in [1,\infty)$, one has
$L^p(\Omega,\lambda)' = L^q(\Omega,\lambda)$ with $q^{-1} + p^{-1} = 1$. In 
particular, this is true with $p =1$.
\end{proposition}
As a consequence, one has
\begin{corollary}\label{cor:Yosida}
Assume that $\CF$ is countably generated
and let $g_n$ be a bounded sequence in $L^\infty(\Omega,\lambda)$. Then there exists
a subsequence $g_{n_k}$ and an element $g \in L^\infty(\Omega,\lambda)$ such that
$\int g_{n_k}(x)f(x)\,\lambda(dx) \to \int g(x) f(x)\,\lambda(dx)$ for every $f \in L^1(\Omega,\lambda)$.
\end{corollary}
\begin{proof}
Since $\CF$ is countably generated, $L^1(\Omega,\lambda)$ is separable and therefore
contains a countable dense subset $\{f_m\}$. Since the $g_n$ are uniformly bounded,
a diagonal argument allows to exhibit
an element $g \in L^1(\Omega,\lambda)'$ and a subsequence $n_k$ such that
$\int g_{n_k}(x)f_m(x)\,\lambda(dx) \to \scal{f_m,g}$ for every $m$. The claim follows
from the density of the set $\{f_m\}$ and the previous proposition.
\end{proof}

Note also that one has

\begin{lemma}\label{lem:refmeas}
Let $P$ be a strong Feller Markov kernel on a Polish space $\CX$. Then there exists
a probability measure $\lambda$ on $\CX$ such that $P(x,\cdot\,)$ is absolutely
continuous with respect to $\lambda$ for every $x \in \CX$.
\end{lemma}

\begin{proof}
Let $\{x_n\}$ be a countable dense subset of $\CX$ and define a probability measure $\lambda$ by 
$\lambda(A) = \sum_{n=1}^\infty 2^{-n} P(x_n,A)$.
Let $x\in\CX$ be arbitrary and assume by contradiction that $P(x,\cdot\,)$ is not
absolutely continuous with respect to $\lambda$. This implies that there exists a set
$A$ with $\lambda(A) = 0$ but $P(x,A) \neq 0$. Set $f = \chi_A$ and consider $Pf$.
One one hand, $\bigl(Pf\bigr)(x) = P(x,A) > 0$. On the other hand, $\bigl(Pf\bigr)(x_n) = 0$
for every $n$. Since $P$ is strong Feller, $Pf$ must be continuous, thus leading to a contradiction.
\end{proof}

Finally, to complete the preliminaries, set $B = \{g \in \CB_b(\CX)\,|\,\sup_x |g(x)| \le 1\}$ 
the unit ball in the space of bounded measurable functions, and note that one has the
following alternative formulation of the ultra Feller property:

\begin{lemma}\label{lem:equicont}
A Markov kernel $P$ on a Polish space $\CX$ is ultra Feller if and only if the set of functions
$\{Pg\,|\, g \in B\}$ is equicontinuous.
\end{lemma}

\begin{proof}
This is an immediate consequence of the fact that one has the characterisation
$\|P(x,\cdot\,) - P(y,\cdot\,)\|_{\mathrm{TV}} = \sup_{g\in B} |Pg(x) - Pg(y)|$, see for example 
\cite[Example~1.17]{Vil}.
\end{proof}

We have now all the ingredients necessary for the proof of the result announced earlier.

\begin{theorem}\label{theo:SF}
Let $\CX$ be a Polish space and let $P$ and $Q$ be two strong Feller
Markov kernels on $\CX$. Then the Markov kernel $PQ$ is ultra Feller.
\end{theorem}

\begin{proof}
Applying Lemma~\ref{lem:refmeas} to $Q$, we see that there exists a reference measure $\lambda$
such that $Q(y,dz) = k(x,z)\,\lambda(dz)$.

Suppose by contradiction that $R = PQ$ is not ultra Feller. Therefore, by Lemma~\ref{lem:equicont} there
exists an element $x \in \CX$, a sequence $g_n \in B$, a sequence $x_n$ converging to $x$, and
a value $\delta > 0$ such that 
\begin{equ}[e:step]
Rg_n(x_n) - Rg_n(x) > \delta\;,
\end{equ}
for every $n$. 
Interpreting the $g_n$'s as elements of
$L^\infty(\CX,\lambda)$, it follows from Corollary~\ref{cor:Yosida} that, 
extracting a subsequence if necessary, we can assume that there exists an element
$g \in L^\infty(\CX,\lambda)$ such that
\begin{equ}
\lim_{n \to \infty} Qg_n(y) = \lim_{n \to \infty} \int k(y,z) g_n(z)\,\lambda(dz) 
=  \int k(y,z) g(z)\,\lambda(dz) = Qg(y)
\end{equ}
for every $y \in\CX$. (This is because $k(y,\cdot\,) \in L^1(\CX,\lambda)$.)
Let us define the shorthands $f_n = Qg_n$, $f = Qg$, and $h_n = \sup_{m\ge n} |f_m - f|$.

Since $f_n \to f$ pointwise, it follows from Lebesgue's dominated convergence theorem
that $Pf_n(x) \to Pf(x)$. The same argument shows that $Ph_n(y) \to 0$ for every $y \in \CX$. Since 
furthermore the $h_n$'s are positive decreasing functions, one has
\begin{equ}
\lim_{n \to \infty}Ph_n(x_n) \le \lim_{n \to \infty}Ph_m (x_n) = Ph_m(x)\;,
\end{equ}
which is valid for every $m$, thus showing that $\lim_{n \to \infty}Ph_n(x_n) = 0$.
This implies that 
\begin{equs}
\lim_{n \to \infty} Pf_n(x_n) - Pf(x) &\le \lim_{n \to \infty} |Pf_n(x_n) - Pf(x_n)| + \lim_{n \to \infty} |Pf(x_n) - Pf(x)| \\
&\le \lim_{n \to \infty} Ph_n(x_n) + 0 = 0\;,
\end{equs}
thus creating the required contradiction with \eref{e:step}.
\end{proof}

\begin{example}
Let us conclude this section with an example of a strong Feller Markov kernels which is not ultra Feller.
Take $\CX = [0,1]$ and define $P$ by
\begin{equ}
P(x, dy) = \left\{\begin{array}{cl} dy & \text{if $x = 0$,} \\ 
c(x)\bigl(1 + \sin(y/x)\bigr)\, dy& \text{otherwise.} \end{array}\right.
\end{equ}
Here, the function $c$ is chosen in such a way that $P(x,\cdot\,)$ is a probability measure.
It is obvious that, for any $f \in \CB_b(\CX)$, $P f$ is continuous (even $\CC^\infty$) 
outside of $x=0$. It follows furthermore from the Riemann-Lebesgue lemma that $P f$ is continuous
at $x = 0$. However, the map $x \mapsto P(x,\cdot\,)$ is discontinuous at $0$ in the total variation
topology (one has $\lim_{x \to 0} \|P(x,\cdot\,) - P(0,\cdot\,)\|_\TV = {2\over \pi}$), which shows that $P$ is not ultra Feller. 
\end{example}

\begin{remark}
Since, as seen in the previous example, there are strong Feller Markov kernels that are not ultra Feller, 
Theorem~\ref{theo:SF} fails in general if one of the two kernels is only Feller (take the identity). 
\end{remark}

\section[Some Gaussian measure theory]{Appendix B: Some Gaussian measure theory}
\label{sec:Gaussian}

This section is devoted to a short summary of the theory of Gaussian measures and in particular on their conditioning.
Denote by $X$ some separable Fr\'echet space and assume that we are given a splitting $X = X_1 \oplus X_2$. 
This means that the $X_i$ are subspaces of $X$ and every element of $X$ can be written uniquely as $x = x_1 + x_2$
with $x_i \in X_i$ and the projection maps $\Pi_i \colon x \mapsto x_i$ are continuous.

Assume that we are given a Gaussian probability measure $\P$ on $X$, with covariance operator $Q$. That is $Q \colon X^*  \to X$ is a continuous bilinear map such that $\scal{Qf,g} = \int f(x)g(x)\,\P(dx)$ for every $f$ and $g$ in $X^*$. 
(Such a map exists because $\P$ is automatically a Radon measure in our case.) Here, we
used the notation $\scal{f,x}$ for the pairing between $X^*$ and $X$.
We denote by $\CH$ the reproducing kernel Hilbert space of $\P$. The space $\CH$ can be constructed
as the closure of the image of the canonical map $\iota \colon X^* \to L^2(X,\P)$ given by
$\bigl(\iota h\bigr)(w) = h(w)$, so that $\CH$ is the space of $\P$-measurable linear functionals on $X$. 
If we assume that the support of $\P$ is all of $X$ (replace $X$ by the support of $\P$ otherwise), then 
this map is an injection, so that we can identify $X^*$ with a subspace of $\CH$.
Any given $h \in \CH$ can then be identified with the (unique)
element $h_*$ in $X$ such that $\scal{Qg, h} = g(h_*)$ for every $g \in X^*$.
With this notation, the scalar product on $\CH$ is given by
$\scal{\iota h, w} = h(w)$, or equivalently by $\scal{\iota h, \iota g} = \scal{h, Qg} = \scal{Qh, g}$. We will from now on use 
these identifications, so that one has
\begin{equ}[e:ident]
X^* \subset \CH \subset X\;,
\end{equ}
and, with respect to the norm on $\CH$, the map $Q$ is an isometry between $X^*$ and its image.
For elements $x$ in the image of $Q$ (which can be identified with a dense subset of $\CH$), one has $\|x\|^2 = \scal{x, Q^{-1} x}$.

Given projections $\Pi_i \colon X \to X_i$ as above, the reproducing kernel Hilbert spaces $\CH_i^p$ of the projected
measures $\P \circ \Pi_i^{-1}$ are given by $\CH_i^p = \Pi_i \CH \subset X_i$, and their covariance operators
$Q_i$ are given by $Q_i = \Pi_i Q \Pi_i^* \colon X_i^* \to X_i$. The norm on $\CH_i^p$ is given by
\begin{equ}
\|x\|_{i,p}^2 = \inf \{ \|y\|^2 \,:\, x = \Pi_i y\;,\; y \in \CH\} = \scal{x, Q_i^{-1} x}\;,
\end{equ}
where the last equality is valid for $x$ belonging to the image of $Q_i$. It is noteworthy that even though the spaces
$\CH_i^p$ are not subspaces of $\CH$ in general, there is a natural isomorphism between $\CH_i^p$ some closed
subspace of $\CH$ in the following way. For $x$ in the image of $Q_i$, define $U_i x = Q \Pi_i^* Q_i^{-1} x \in X$.
One has
\begin{lemma}\label{lem:isom}
For every $x$ in the image of $Q_i$, one has $U_i x \in \CH$. Furthermore, the map $U_i$ extends to an isometry
between $\CH_i^p$ and $U_i \CH_i^p \subset \CH$.
\end{lemma}

\begin{proof}
Since $U_i x$ belongs to the image of $Q$ by construction, one has $\|U_i x\|^2 = \scal{U_i x, Q^{-1} U_i x} = \scal{Q \Pi_i^* Q_i^{-1} x, \Pi_i^* Q_i^{-1} x} = \scal{x, Q_i^{-1} x} = \|x\|_{i,p}^2$.  The claim follows from the fact that the image of $Q_i$ is
dense in $\CH_i^p$.
\end{proof}

We denote by $\hat \CH_i^p$ the images of $\CH_i^p$ under $U_i$. Via the identification \eref{e:ident}, 
it follows that $\hat \CH_i^p$ is actually nothing but the closure in $\CH$ of the image of $X_i^*$ under $\Pi_i^*$.
Denoting by $\hat \Pi_i \colon \CH \to \CH$ the orthogonal projection (in $\CH$) onto $\hat \CH_i^p$, it is a straightforward
calculation to see that one has the identity $\hat \Pi_i x = U_i \Pi_i x$. On the other hand, it follows from the definition of $Q_i$
that $\Pi_i U_i x = x$, so that $\Pi_i \colon \hat \CH_i^p \to \CH_i^p$ is the inverse of the isomorphism $U_i$.

We can also define subspaces $\CH_i^c$ of $\CH$ by
\begin{equ}[e:carCHc]
\CH_i^c = \overline{\CH \cap X_i} =  \overline{\CH \cap \CH_i^p} \;,
\end{equ}
where we used the identification \eref{e:ident} and the embedding $X_i \subset X$. The closures are taken
with respect to the topology of $\CH$. The spaces $\CH_i^c$ are
again Hilbert spaces (they inherit their structure from $\CH$, not from $\CH_i^p$!) and they therefore define 
Gaussian measures $\P_i$ on $X_i$. Note that for $x \in \CH_i^c \cap \CH_i^p$, one has $\|x\| \ge \|x\|_{i,p}$, so that the
inclusion $\CH_i^c \subset \CH_i^p$ holds. One has
\begin{lemma}\label{lem:ortho}
One has $\CH_1^c = \bigl(\hat \CH_2^p\bigr)^\perp$ and vice-versa.
\end{lemma}
\begin{proof}
It is an immediate consequence of the facts that $X = X_1 \oplus X_2$, that $\hat \CH_1^p$ is the closure of the image
of $\Pi_1^*$ and that, via the identification \eref{e:ident} the scalar product in $\CH$ is an extension of the duality pairing
between $X$ and $X^*$.
\end{proof}

We now define a (continuous) operator $P \colon \CH_1^p \to \CH_2^p$ by $Px = \Pi_2 U_1 x$. 
It follows from the previous remarks that $P$ is unitarily equivalent to the orthogonal projection (in $\CH$) from 
$\hat \CH_1^p$ to $\hat \CH_2^p$. Furthermore, one has $P x = U_1 x - x$, so that
\begin{equ}[e:boundP]
\|Px\|_\CH \le \|x\|_{\CH_1^p} + \|x\|\;,
\end{equ}
which, combined with \eref{e:carCHc}, shows that $P$ can be extended to a bounded operator
from $\CH_1^c$ to $\CH_2^c$.

A standard result in Gaussian measure theory states that $P$ can be extended 
to a $\bigl(\P\circ \Pi_1^{-1}\bigr)$-measurable linear operator $\hat P \colon X_1 \to X_2$. With these notations at hand,
the main statement of this section is given by:

\begin{proposition}\label{prop:conditional}
The measure $\P$ admits the disintegration
\begin{equ}[e:dis]
\int \phi(x) \, \P(dx) = \int_{X_1} \int_{X_2} \phi(x + \hat Px + y) \P_2(dy) \, \bigl(\P\circ \Pi_1^{-1}\bigr)(dx)\;.
\end{equ}
\end{proposition}

\begin{proof}
Denote by $\nu$ the measure on the right hand side.
Since $\nu$ is the image of the Gaussian measure $\mu = \bigl(\P\circ \Pi_1^{-1}\bigr) \otimes \P_2$ under the
$\mu$-measurable linear operator $A \colon (x,y) \mapsto x + \hat P x + y$, it follows from
\cite[Theorem~3.10]{Bog98Gauss} that $\nu$ is again a Gaussian measure. The claim then follows if we can show that the reproducing kernel
Hilbert space of $\nu$ is equal to $\CH$. Since the reproducing kernel space $\CH(\mu)$ of $\mu$ is canonically isomorphic to
$\CH(\mu) = \CH_1^p \oplus \CH_2^c \subset \CH \oplus \CH$, this is equivalent to the fact that the operator $x \mapsto x + P x = x + \Pi_2 U_1 x$
from $\CH_1^p$ to $\CH$ is an isometry between $\CH_1^p$ and $\bigl(\CH_2^c\bigr)^\perp$. On the other hand, we know from Lemma~\ref{lem:ortho}
that $\bigl(\CH_2^c\bigr)^\perp = \hat \CH_1^p$ and we know from Lemma~\ref{lem:isom}
that $U_1$ is an isomorphism between $\CH_1^p$ and $\hat \CH_1^p$. 
Finally, it follows from the definitions that $\Pi_1 U_1 x = Q_1 Q_1^{-1} x = x$ for every $x\in \CH_1^p$, 
so that
one has $x + \Pi_2 U_1 x  = \bigl(\Pi_1 + \Pi_2\bigr) U_1 x = U_1 x$, which completes the proof.
\end{proof}


\makeatletter
\let\mean\@old@mean
\let\d\@old@d
\let\cM\@old@cM
\let\cP\@old@cP
\let\cF\@old@cF
\let\P\@old@P
\let\E\@old@E
\let\R\@old@R
\let\Z\@old@Z
\let\N\@old@N
\let\CC\@old@CC
\let\CH\@old@CH
\let\CB\@old@CB
\let\CF\@old@CF
\let\CM\@old@CM
\let\CP\@old@CP
\let\CX\@old@CX
\let\CY\@old@CY
\let\CQ\@old@CQ
\let\CS\@old@CS
\let\CW\@old@CW
\let\TV\@old@TV
\let\eref\@old@eref
\let\scal\@old@scal
\makeatother

\end{document}